\documentclass{article}
\usepackage[latin1]{inputenc}
\usepackage[T1]{fontenc}
\usepackage[english]{babel}
\usepackage{amsmath}
\usepackage{amssymb,amsfonts,textcomp}
\usepackage{color}
\usepackage{array}
\usepackage{xtabular}
\usepackage{hhline}
\usepackage{hyperref}

\date{}

\begin{document}
\title{India's Contribution to Arab mathematics\footnote{Translated by Dileep Karanth, Department of Physics, University of Wisconsin-Parkside, Kenosha.}}
\author{Khalil Jaouiche}
%\author{Translated by Dileep Karanth}
%\affiliation{ Department of Physics, University of Wisconsin-Parkside, Kenosha, WI 53141}
%\date{\today}
\maketitle

\long\def\symbolfootnote[#1]#2{\begingroup%
\def\thefootnote{\fnsymbol{footnote}}\footnote[#1]{#2}\endgroup} 

\begin{abstract}
This paper was presented at the Colloque de Saint-Denis de la Reunion, November 3-7, 1997, and published  in the book \textit{L'Océan Indien au carrefour des mathématiques arabes, chinoises, européennes et indiennes} (pp. 211-223), Tournès, Dominique (ed.), Saint-Denis, Réunion, I.U.F.M. de La Réunion, 1998. 

The article is an introductory survey of the influence of India on the mathematics of the Islamic world. 

\end{abstract}

The study of the connections between Arab and Indian mathematics is beset by a number of difficulties. The first kind of difficulty is due to the paucity of documentation. There is no work or study from the Arab Middle Ages that deals with this question. The only general book written about India by an author of the Arab-Muslim world, that by Al-B\={\i}r\={u}n\={\i}, written in 1030,\footnote[1]{\textit{Kit\={a}b ta\d{h}q\={\i}q m\={a} li-l-Hind min maq\={u}latin fi l-'aql aw marz\={u}latin} (Enquiry into what is said in India which is conformity with or in opposition to reason)} dedicates several pages to Indian astronomy, but does not treat either  arithmetic or algebra.  The second kind of difficulty has to do with the extreme complexity of the relations between the Arab and Middle Eastern world on the one hand, and India and China on the other. The studies in this colloquium have brought to light the existence of an Indo-Chinese cultural ensemble within which it is possible to detect similarities in methods of calculation, but not to establish a definitive chronology between these methods. To limit the study to the connections between Arab and Indian mathematics amounts to excluding quite a bit that these schools of mathematics owe directly or indirectly to China.

If we add to these difficulties the existence of an influence of Greek mathematics since the third century BC, on the Near East without doubt, but also probably on China and India, the connections between Arab and Indian mathematics seem like an island that can be separated only artificially from the multiple and complex relationships going from the Mediterranean to the China Sea, with the Indian Ocean as a cross-roads, as the theme of this colloquium aptly reminds us.

However, we will proceed to make such an artificial sundering in the lines that follow. Only focussed studies can progressively disentangle the web of these complex connections. We are going to devote ourselves to a study of some chapters of Arab mathematics where the influence of India is explicitly or implicitly attested.

\section{}
That India has played a preponderant role in the formation of Islamic culture, since the days of its appearance in the VIIIth century of the Christian era, can be seen in the fact that the first great classical work of Arabic literature -- a book of fables titled \textit{Kal\={\i}l\={a} wa Dimn\={a}} -- is the translation from Pahlavi of an Indian work whose protagonists are two animals of which one represents hypocrisy and greed and the other honesty and sincerity, with the latter finally prevailing over the former. 

However, it was through the medium of astronomy that Indian science first penetrated the Arab world. In 772, following the visit of a great Indian astronomer to the court of Bagdad, the caliph al-Mansour (754-775) charged al-Faz\={a}r\={\i} with the translation into Arabic of the astronomical tables of the \textit{Siddh\={a}nta}. This translation, which has come down us only in fragments,\footnote[2]{F. Sezgin, Geschichte des Arabischen Schriftums, vol. VI, Leiden, 1978, p. 123)} is known in Arabic by the name 
\textit{as-Sindhind al-Kab\={\i}r}. While al-Kab\={\i}r in Arabic means ``big'' or ``great'', the word \textit{sindhind} in Indian\symbolfootnote[1]{Translator's Note: The word ``Indian'' here  probably refers to the Sanskrit language.} [sic], according to ibn al-Qifti's chronicle, means \textit{ad-dahr ad-d\={a}hir}\footnote[3]{Ibn al-Qifti, \textit{Ta'ri\underline{kh}}, Lippert, Leipzig, 1903, p. 270.}, expression equivalent to our ``centuries of centuries''. According to Pingree, who has studied these fragments\footnote[4]{Sezgin, \textit{ibid.}}, they contain problems dealt with by Greek astronomy. 

Indian astronomy would again make its appearance in the Middle east, but this time through the pen of the first great mathematician of Islam: Mohammad ibn Mussa Al-Khwarizm\={\i}. Hailing from Khwarizm, the region situated south of the Aral Sea, as his name indicates, he was one of the most active and celebrated scholars of the ``House of Wisdom'' at Bagdad. This institute, founded by al-Ma'mun, functioned as a veritable centre of research where the most brilliant of the translators and scholars of Baghdad came together. It was there that Al-Khwarizm\={\i} edited the two versions of the \textit{zij as-Sindhind }(astronomical tables) that Ibn al-Nad\={\i}m\footnote[5]{Al-Fihrist, Fl\"ugel, p. 274} attributes to him. Of it today we only have the Latin translation of Adelard\footnote[6]{Spelt as ``Ath\'{e}lard'' in German works.} of Bath, which was probably based on the second version, revised by Abul-Q\=asim al-Maghr\={\i}ti (d. 1007 CE)\footnote[7]{H. Suter, \textit{Die Astronomischen Tafeln des Muhammad Ibn Musa al-Khwarizmi}, edited by the Institut f\"ur Geschichte der Arabisch-Islamischen Wissenschaften, Frankfurt, 1986, in Suter, vol. 1, p. 473 ff.}. The original version of this \textit{zij} being lost, and since we only possess fragments from al-Faz\={a}r\={\i}'s tables, it is difficult for us to determine what Al-Khwarizm\={\i}'s true sources were. It seems that the Indian sources were the \textit{Br\={a}hma-Sphu\d{t}a-Siddh\={a}nta} of Brahmagupta (1st half of the VIIth Century CE) and the \textit{S\=urya-Siddh\={a}nta}, on which were based the \textit{zij as-Shah}, translated from Persian into Arabic in the VIIIth CE. However, as Al-Khwarizm\={\i}'s text shows, it is necessary to add to these sources the influence of Ptolemy's Almageste as well as the results of the observations made by the astronomers of the ``House of Wisdom". On the mathematical plane, the use of trigonometric lines: sine, cosine, tangent and cotangent is current in the text\footnote[8]{See, for a discussion of this question, Suter, \textit{ibid.}, as well as F. Sezgin, \textit{op. cit.}, vol., VI, p. 141-143.}.

\section{}
It is also to Al-Khwarizm\={\i} that we owe the second work in which Indian influence is manifest. This is his book on arithmetic, which he explicitly borrows from India. The exact title of the book remains unknown, both in the Arabic original, now lost, as well as in the translations or Latin versions of the XIIth century which have come down to us, and of which Allard has made a remarkable critical edition, accompanied by the French translation\footnote[9]{A. Allard, \textit{Muhammad Ibn Musa al-Khwarizmi, Le calcul indien}, Librairie scientifique et technique Albert Blanchard, Paris, 1992.}.

The importance of this work for the history of mathematics, both in the West and in the East, cannot be underestimated. It was this book that gave decimal numeration the eminent place that it has since occupied in the Mediterranean basin. 

For the story, let us recall that it begins with two words, of which one would go on to enjoy an exceptional fortune in the history of mathematics, and today, in computer science: ``Dixit Algorizmi'' (Al-Khwarizm\={\i} has said). This is the origin of our word ``algorithm'' used today to designate any operational recipe in a given calculus\footnote[10]{\textit{Le calcul indien},\textit{ibid.}, p. 1. It is by this name we will refer to the edition and translation by M. Allard, in the notes that follow.}.

The work has a double objective. First, to explain the system of numeration of the Indians ``who use 9 letters for all their numbers thanks to a method which is peculiar to them''\footnote[11]{\textit{Ibid.}, p. 1-6.}. This is an explanation of the decimal system, where, as everyone knows, each of the nine digits that compose it, designated by a ``letter'', has a particular value according to the position it occupies: that of units, tens, hundreds, thousands, tens of thousands, etc. 

It is in this exposition of the properties of the decimal system that Al-Khwarizm\={\i}  describes the origin of what would become our ``zero''. One could say that neither the Arabs nor the Indians nor the Chinese ``invented'' the zero. As we now conceive it, it is the fruit of a long evolution which culminated in the XIXth century with the publication in 1889, by Giuseppe Peano, of his \textit{Arithmetices principia nova methodo exposita}, which contained his famous axioms. The first of them affirms: ``Zero is a number''\footnote[12]{Cf. Carl Boyer, \textit{A History of Mathematics}, Princeton University Press, Wiley \& Sons, 1968, p. 645.}.

Such a conception would have profoundly shocked Al-Khwarizm\={\i} as much as his Indian and Chinese predecessors. Not only would it have been strange to them, but they did not have any word for what \textit{we} now call ``zero''. This why, we think, the use of this word by Allard throughout his remarkable translation is unfortunate. The word ``zero'', for a modern reader, carries the connotation of a number that possesses certain properties defined by the rules to which it is subject in the calculus. It is neither positive, nor negative; it is smaller than all the positive reals, and larger than all the negative real numbers, etc. So many notions foreign to the Arabs, the Indians and doubtless to the Chinese. The zero is is for us a mathematical object that really exists, whereas as we shall see immediately, in Al-Khwarizm\={\i} 's text, it is a pure non-object, a ``void''. This is the reason why we have replaced the word ``zero'' by the word ``ring''\symbolfootnote[1]{Translator's Note: The word ``ring" is an attempt to translate the French word ``rond".}, literal translation of the Latin ``circulus''\footnote[13]{``Circulus'' is without doubt the Latin translation of the Arabic \textit{d\=a'ira}. In Arabic, this word stands for ``circle'' as well as for ``ring''. But Al-Khwarizm\={\i} does not refer to the geometric figure in so many words, but as he says explicitly a few lines later, to a sign that ``resembles the letter O''.}, in the excerpt from Allard's translation, that we will quote shortly.

Right in the beginning of his work, Al-Khwarizm\={\i} recalls justifiably ``that in no position there is more than 9 nor less than 1, except for the case when there is a ring, which is nothing\footnote[14]{\textit{Le calcul indien}, p. 6; emphasis ours.}. The decimal system can only comprise nine digits. Moreover, as we have said before, Al-Khwarizm\={\i} clearly says that the value of a digit depends on the place it occupies in the number, counting positions from the right (units, tens, hundreds ...). The first difficulty which arises in this system is the writing of the number 10 -- and then the tens -- which cannot be represented by a single ``letter'' (digit). It can only be represented by a 1 placed in the tens' position. But how could one specify this position, in a writing that did not trace columns to indicate position? Al-Khwarizm\={\i} tells us the answer\footnote[15]{\textit{Ibid.}, p.3.}: ``... a representation of tens was necessary for them (the Indians) because it was similar to the representation of one, so that one could recognize by it that by it 10 was meant. So they placed in front of it a position, and placed in it a little ring \textit{resembling the letter O to know by it that the position of units was empty}\footnote[16]{Emphasis ours.}, that it contained no number except for the little ring of which we have said it occupied the place, and the number which occupied the next place was a ten...'' The same problem arises evidently when we try to write a number, one of whose places ``does not contain anything.'' In that case, ``you place a ring so that the position is not blank, but is occupied by a ring, so as not to the reduce the places, and one does does not think that the second is the first, and in that way, be mistaken about the number''\footnote[17]{\textit{Ibid.}, p.7.}. Thus, for example, one wishes to write out in digits the number two hundred and four, it is theoretically enough to write a 4 and, a little further, on the left, a 2. But then one could possibly read it as twenty-four and not two hundred and four. So, \textit{as a precaution}, one places a ``ring'' between the 4 and the 2 in order to indicate clearly that the latter occupies the hundreds place.

	To be sure, the word ``ciffre'', from which the word ``zero'' is derived etymologically, as well as the French word ``chiffre'' (meaning digit, or number), is not totally absent in {\textit{Le calcul indien}. It figures in it at least once with its \textit{two} possible symbols. ``They [the Indians] also use'', says Al-Khwarizm\={\i}, "the sifr in this manner -- O and $\tau$''\footnote[18]{\textit{Op. cit.}; Allard had translated by ``zero'', as for ``circulus''.}. But the word \textit{sifr}, translated into Latin as ``ciffre'' has the meaning ``empty'' in Arabic, as noted by Ifrah\footnote[19]{\textit{Histoire universelle des chiffres}, Seghers, Paris, 1981, (first edition), p. 509.}, in specifying that it was the Arabic translation of the Sanskrit ``\'s\=unya'', which has precisely the same meaning. And Ifrah, citing Taton\footnote[20]{R.Taton, \textit{Histoire du calcul}, Coll. ``Que sais-je ?'', P.U.F., Paris, 1969. }, recalls that ``in the XIIIth century, in France, popular language qualified a man devoid of culture as \textit{Cyfre d'angorisme} or even as \textit{Cifre en algorisme}''\footnote[21]{Ifrah, \textit{op. cit.}, p. 512}.
	
	Is it surprising that the Indians had intuitively deduced that the result of a multiplication of a ``ring that stands for nothing'' (\textit{circulus nihil significans})\footnote[22]{\textit{Le calcul indien}, p. 29.} with a number, is nothing? Al-Khwarizm\={\i}\footnote[23]{\textit{Ibid.}, p.10.} writes, ``any ring multiplied by any number whatsoever, is nothing, that is, no number results from it, that anything that is multiplied by a ring is also nothing \ldots''. To write, as is habitually done, that the Indians, the Arabs and the Chinese knew that the product of an arbitrary number by zero is equal to zero\footnote[24]{Cf. for example, Youschkevitch, \textit{Les math\'ematiques arabes}, Vrin, Paris, 1976, p. 18.} is to obscure the picture or the conception that these peoples had of such an operation, and to substitute in its place a conception which would be established only much later. 
	
	The second aim of the {\textit{Le calcul indien} is to explain the calculation methods of the elementary arithmetic operations: addition, subtraction, multiplication, division, the extraction of square roots, multiplication of decimal and sexagesimal fractions, duplication\footnote[25]{An operation which consists in carrying out a multiplication more easily by writing the factor in the form of a sum of powers of 2, and the multiplying the multiplicand by each of these powers. Cf. M. Caveing, \textit{Essai sur le savoir math\'ematique dans la M\'esopotamie et l'\'Egypte ancienne},Presses Universitaires de Lille, 1994, p. 251-253.} and dimidiation\footnote[26]{A method of division inverse to the multiplication, using successive powers of $\frac{1}{2}$.}. To give the reader an idea of the calculation methods used in the Near East at least from the IXth century onwards -- and much later than in India -- and the form in which they passed into the Occident from the XIIth century onwards, we are going to cite an example given by  Al-Khwarizm\={\i}\footnote[27]{\textit{Le calcul indien, op. cit.}, p. 9-11. The different stages of the calculation being quite difficult to follow in the text, we are going to borrow them from A. Youschkevitch,  \textit{Les math\'ematiques arabes}, p. 19.}. 
\newpage
Consider the multiplication of 2326 by 214.
\begin{enumerate}
\item We write
\begin{center}
{ \ \ \ \ \ \ \ }  2326 \\
214
\end{center}
with the units place of the multiplier (4) placed under the highest ranked place of the multiplicand (2). 
\item We multiply the highest ranked position of the  multiplicand (2) by the multiplier: $ 2 \times 214 = 428 $. 
\item We write the product to the left of the multiplicand in such a way that the 8 in the units place of 428 replaces the 2 which has already been used for multiplication, and we move the multiplier one place to the left. We then have the following two lines:
\begin{center}
{ \ \ \ \ \ }  428326 \\
{ \ \ \  } 214
\end{center}
\item Now it is the number 3 of the multiplicand that is to be multiplied by the multiplier: $ 3 \times 214 = 642 $. 
\item This last product is now written over the multiplicand in such a way that the 2 in the units place of the number 642 is just above the 3 which has just been multiplied by 214; We then have the following lines:
\begin{center}
{ \ \ \  } 642 \\
{ \ \ \ \ \ }  428326 \\
{ \ \ \  } 214
\end{center}
\item The two in the units place of 642 must replace the 3 below it. We add 64 and 28, and we move the multiplier one place to the left, to get:
\begin{center}
{ \ \ \ \ \ }  492226 \\
{ \ \ \  } 214
\end{center}
\item We multiply the 2 which earlier came after the 3, and thus get the following lines
\begin{center}
{ \ \ \  \ } 428 \\
{\ \ \ }  492226 \\
{ \ \ \  \ } 214
\end{center}
\item We perform the same operation as before and get the following lines:
\begin{center}
{ \ \ \ \ \ }  496486 \\
{ \ \ \ \ \ \ \ \ \  } 214
\end{center}
\item Finally we multiply the 6 in the units place in the upper line -- the last digit of the initial multiplicand -- by 214, to get 1284 -- and we perform once more the same operation as before, and we the following two lines:
\begin{center}
{ \ \ \ \ \ }  497764 \\
{ \ \ \ \ \ \ \ \ \  } 214
\end{center}
\item The upper line is the result of the multiplication.
\end{enumerate}

However {\textit{Le calcul indien} is of historical interest not only because it is a complete treatise dealing with calculation methods. There are even some surprising reminiscences from Greek philosophy and geometry to be found in it. Thus, there is a recollection of Aristotelian notions concerning unity as a principle of the science of numbers, different sorts of movements, concepts of diminution and augmentation, and also the \textit{zij} of Ptolemy, in a jumbled form in a passage dealing with the multiplication of sexagesimal fractions\footnote[28]{\textit{Op. cit.}, p. 23-24.}. It is not very likely that these notions had been introduced into to the {\textit{Le calcul indien} by Al-Khwarizm\={\i} himself based on contemporary translations made by translators working at the House of Wisdom (Bait al-\d{H}ikma). These translations would have presented coherent and structured texts, whereas Al-Khwarizm\={\i} juxtaposes notions dealt with by Aristotle in various different works, which have no direct relation to each other. The translation of the work by Ptolemy bears the name \textit{al-Majisti} in the Arabic translations of the time, and not the name \textit{zij}. All this seems to indicate we must be dealing with the relics of Greek science and philosophy, which had passed over to India well before the advent of Islam.

	The same is the case with the mention of evenly even, evenly odd and oddly even, numbers defined by Nicomachus of Gerasa\footnote[29]{I thank Mr. Bernard Vitrac for having brought to my notice the distinction between the three kinds of numbers, as they are defined in \textit{Le calcul indien}, harking back to Nicomachus of Gerasa and not to definitions 8-10 of Book VII of Euclid's Elements, where there is in addition a fourth kind of number: the oddly odd.} (Ist-IInd century AD), to which the author adds some rather trivial explanations. It is hard to see why Al-Khwarizm\={\i} would have introduced these kinds of numbers if he had not found them in the Indian original. But here once again, there arises the question of the relationship between Greek and Indian mathematics.
	
	It will be seen that while the {\textit{Le calcul indien} by Al-Khwarizm\={\i}, whose historical importance is considerable, is an example of the incontestable influence of Indian mathematics on that of the Islamic world, and hence, on that of the occidental Middle Ages, it brings to light the methodological difficulties which arise when one tries to isolate a culture to study its influence on another. In the form that it has acquired, {\textit{Le calcul indien} is a kind of melting pot in which cultural elements from various civilisations -- Babylonian, Egyptian, Greek, 
Indian and Chinese -- have all fused together. It is difficult to establish a strictly linear chronology between them. 

\section{}
	The {\textit{Le calcul indien} is an important example which illustrates the close relations that must have existed between the Arab-Islamic world and India, but it is by no means the only one. We will have occasion again to find Indian mathematics at a later age -- the age of Al-Karaj\={\i} (died circa 1029-1030 AD). Al-Karaj\={\i} paid much attention to addition, subtraction and to the extraction of the square roots of irrational numerical polynomials, matters in which the Indian mathematician Bh\=askara II (1114-1185) also distinguished himself\footnote[30]{There also was a Bh\=askara I, who lived in the VIth century of the Christian era.}. 
	
	We will first give the example of the extraction of the square root of the sum of a number of irrational roots. Suppose we wish to calculate\footnote[31]{Al-Karaj\={\i}, \textit{al-Badi`}Edited with introduction (in French) and notes by Adel Anbouba, Publications de l'Universit\'e libanaise, section des math\'ematiques, Beirut, 1964, p. 44-45.},:
\begin{center} $ \sqrt{16 + \sqrt{24}  + \sqrt{40}  + \sqrt{48}  + \sqrt{60}  + \sqrt{72}  + \sqrt{120}} $ \end{center}
	This expression consists of seven terms of which the first, 16, which is rational, is the sum of a certain number of squares. The other terms, six in all, are nothing but the double product of the roots taken two at a time. From this Al-Karaj\={\i} concludes the root being sought after can then only have four terms, say:
	\begin{center} $ \sqrt{m}  + \sqrt{n}  + \sqrt{p}  + \sqrt{q}   $ \end{center}
	Al-Karaj\={\i} considered the case when they were arranged in order of increasing magnitude. Then: ``Twenty-four, which is the smallest quantity, is the result of the multiplication of the square of the smallest monomial by the square of the one next to it, four times''. In words: $ 24 = 4 mn$. ``And forty is the result of the multiplication of the square of the smallest monomial  by the square of the third, four times.'' In words: $ 40 = 4 mp$. And, similarly: $ 48 = 4 mq$.``This implies that 24, 40, and 48 are proportional respectively to three monomials of the root\footnote[32]{To be precise, they are proportional to their square.},  with the exception of the first, that is, to $ n, p, q$. Hence, the second of the these three monomials must be a thing.'' Let us set $ n =r $ ($r=$ ``the thing''.) Then $ 4m = 24/r$, and:  $ 40 = (24/r)p$, so that: $ p = 5r/3$. Likewise: $ 48 = (24/r)q$, so: $ q = 2r$. And: $ 60 = 4np$, so that: $ 24/60 = 2/5 = 4mn/4np = m/p$, so: $ m = 2p/5 = 2(5r/3)/5 = (2 \times 5r)/ (3 \times 5) = 2r/3$. But: $ m + n + p + q = 16$, so: $ 2r/3 + 5r/3 + 2r = 16$, so that: $ 16r = 48$, and hence $r =3$. Thus: $ m =2, n =3, p =5, q =6$. The root therefore is:
		\begin{center} $ \sqrt{2}  + \sqrt{3}  + \sqrt{5}  + \sqrt{6}   $ \end{center}
In giving this example, Al-Karaj\={\i} does not refer to Indian mathematics. However, the same polynomial is found in the writings of Bh\={a}skara II (1114-1185), in the \textit{Vijaganita}, or \textit{Bijaganita} according to the spelling adopted by Youschkevitch\footnote[33]{Cf. \textit{Algebra with Arithmetic and Mensuration from the Sanscrit, of Brahmegupta and Bascara}, translated into English by H.T. Colebrook, John Murray, London, 1817, p. 149, and Youschkevitch, \textit{op. cit.}, p. 128. }. The two authors follow different methods for extracting the square root. The method followed by Bh\={a}skara is closer to that of  Al-Karaj\={\i} which we give below. 

	This example illustrates the difficulties entailed in studying the relationships between mathematics in the Islamic world and in India. The fact that Bh\={a}skara had lived more than a century after  Al-Karaj\={\i} would lead one to think that Bh\=askara must have borrowed this example from  Al-Karaj\={\i}. Such a conclusion would be false, however. Indeed, this polynomial is already found in the \textit{Ku\d{t}\d{t}aka}\footnote[34]{Colebrook, \textit{op. cit.}, p.342.} by  Brahmagupta (VIth century BCE). Here we have just one example among many others of the traditionalism which distinguishes all the oriental civilizations. Youschkevitch has also emphasized, regarding the work of Bh\={a}skara II\footnote[35]{\textit{Geschichte der Mathematik im Mittelalter}, Pfalz Verlag, Basel, p. 94-95.}, that the \textit{Siddh\={a}nta-Siroma\d{n}i} (The Crown of the Sciences), of which the \textit{Bija-Ganita} forms a part, has very deep connections with the mathematics of Bh\={a}skara's predecessors. The two mathematicians, Al-Karaj\={\i} and Bh\={a}skara could both have drawn their examples from an earlier tradition hearkening back at least to Brahmagupta. 
	The second example from Al-Karaj\={\i}'s writings that we will cite -- that of the extraction of the square root of an \textit{algebraic polynomial} -- leaves no room for any doubt as to its origin. It is Al-Karaj\={\i} himself who informs us that in his calculations he has followed the method used in the ``Indian reckoning'' (\textit{his\=ab al-Hind}) to extract square roots of ''known quantities'', that if, of numerical polynomials\footnote[36]{Al-Karaj\={\i}, \textit{al-Badi`}, \textit{op. cit.}, text, p.51}. In this example, we have transcribed the calculations in the same terms as used by the Arabs, who obviously were not familiar with our modern symbolism, and used their own terms for degrees of the unknown quantity, namely ``the root'' (denoted by \textit{r}), the dynamis  (denoted by \textit{d}) and the cube (denoted by \textit{c}). The reader can replace $r$ by $x$, $d$ by $x^2$, and $c$ by $x^3$, and easily follow the steps of the calculation.
	Suppose it is desired to extract the square root of:
\begin{center}  $4dcc + 12ddc + 9cc + 20dc + 42dd + 18c +25d + 30r +9 $ \end{center}
	Let us denote the given polynomial, which Al-Karaj\={\i} calls ``the aggregate'' (\textit{al-jumla}), by $P(r)$, and let us denote by $R_j(r)$ the successive remainders of the subtractions which will be carried out.
	\begin{enumerate}
\item Let us take the root or ``the quantity (\textit{al-miqdar}) whose rank is the highest''. We have: $ \sqrt{4dcc} = 2dd $ and $ (2dd) \times (2dd) = 4dcc$. 
\item Let $R_1(r) = P(r) - 4dcc = 12ddc + 9cc + \dots + 9$. 
\item ``You seek the greatest possible quantity which is the closest to 2dd.'' Let A be this this quantity. Then it is necessary that: a) $2A(2dd) $ has a degree equal to the highest degree of $R_1(r)$, that is the degree $ddc$; b) $ A \times A $ must be subtracted from $ R_1(r)$. ``Then you find that'': $ A = 3c$. It then follows that: $2A(2dd) = 12ddc$ and $  A \times A = 9cc$. Then we need to subtract the sum of these two quantities from $R_1(r)$. Then: $ R_2(r) = R_1(r) - (12ddc + 9cc) = 20dc + 42dd + \dots + 9.$ The first two terms of the root are thus:
\begin{center} \framebox{ $2dd + 3c $} \end{center}
\item ``Then you seek the greatest quantity'' B such that: a) $ 2B (2dd + 3c) $ is of a degree equal to the highest degree of $R_2(r)$, that is to the degree of the dynamo-cube; b) $B \times B $ should be subtracted from $R_2(r)$. ``You will then find" that: $R_3(r) = R_2(r) - (20dc + 30dd + 25d) = 12dd + 18c + 30r + 9$. The first three terms of the root, then, are:
\begin{center} \framebox{ $2dd + 3c + 5r$} \end{center} 
\item ``Next seek a quantity'' C, such that: a) $2C (2dd+3c+5r) $ is of a degree equal to the highest degree of $R_3(r)$. Since $R_3(r)$ is already of the degree of the dynamo-dynamis, C can be only a number. ``Thus you will find'' that: C = 3. Then, $ 2C(2dd + 3c + 5r) = 12dd + 18c + 30r$ and $ C \times C = 9$. It is then necessary to subtract the sum of these two quantities from $R_3(r)$. Thus: $R_4(r) = R_3(r) - (12dd + 18c + 30r + 9)$. ``And nothing remains''. The roots of $P(r)$ thus are:
 \begin{center} \framebox{ $2dd + 3c + 5r + 3 $} \end{center}
\end{enumerate}
		The method of calculation used in the preceding example for extracting the square root of this algebraic polynomial is, as we have already said, explicitly borrowed by Al-Karaj\={\i}  from the Indian method of extracting the square root of a \textit{numerical} polynomial. The latter is described by Bh\=askara II in his $ Vija-Ganita$\footnote[37]{Colebrook, \textit{op. cit.}, p. 150-151.} in the chapter dedicated to irrational numbers. 

\section{}
	As-Samaw'al (died ``young'' in 1174), an ardent disciple of Al-Karaj\={\i}, numerous passages of whose writings he quotes in his book \textit{al-B\={a}hir} (The Dazzling), attributes the method of division of two algebraic polynomials to the Indians\footnote[38]{As-Samaw'al, \textit{al-B\=ahir}, ed., French translation and commentary by Roshdi Rashed and Ahmad Salah, publications of the Syrian Ministry for Higher Education, Damascus, 1972, p. 45 ff.}. Of the two examples of polynomial division given by As-Samaw'al, we have chosen the one which contains ``deficient'' quantities (\textit{n\=aqis\=at})\footnote[39]{An expression generally translated incorrectly as ``negative'' quantities. We will return briefly to this question later. We have, however, indicated these quantities by the usual (-) sign, for convenience.}. In the following table, we have labelled the columns with the first letter of quantities listed, as in  Al-Karaj\={\i}'s table. Thus \textit{r }= root,\textit{ d }= dynamis, \textit{c} = cube. Above these appellations we we have listed their modern equivalent, so that the reader can follow the operations more easily. The  quotient, dividend and successive remainder rows contain only the ``coefficients'' of the different terms whose ``degree'' is indicated on the top of the columns. Except for the arrangement of the different steps of the division, the reader can verify that the ``Indian method'' is exactly the same as the method we employ today. 

%\newpage 
\vspace*{0.5in} 

%\begin{flushright}
\tablehead{}
\hspace*{0.90in} 
\begin{xtabular}{|m{0.3in}|m{0.3in}|m{0.3in}|m{0.3in}|m{0.3in}|m{0.3in}|m{0.3in}|m{0.3in}|m{0.3in}|}
\hline
$6x^8$  &
$28x^7$  &
$6x^6$  &
$-80x^5$  &
$38x^4$ &
$92x^3$ &
$-200x^2$  &
$20x$   &
~ \\
$dcc$ &
$ddc$ &
 $cc$ &
$dc$ &
$dd$ &
$c$ &
 $d$ &
$r$ &
units \\
\hline
\end{xtabular}
%\end{flushright}

%\vspace*{1mm}

%\begin{flushright}
\tablehead{}
\begin{xtabular}{|m{1.0in}|m{0.3in}|m{0.3in}|m{0.3in}|m{0.3in}|m{0.3in}|m{0.3in}|m{0.3in}|m{0.3in}|m{0.3in}|}
\hline 
quotient &
~ &
~  &
~  &
~  &
~  &
~  &
~  &
~  &
~  \\\hline
dividend &
6 &
28  &
6  &
-80  &
38  &
92  &
-200  &
20  &
~  \\\hline
divisor &
~ &
~  &
~  &
2  &
8 &
0  &
-20  &
~  &
~  \\\hline
\end{xtabular}
%\end{flushright}

%\vspace*{.1mm}

%\begin{flushright}
\tablehead{}
\begin{xtabular}{|m{1.0in}|m{0.3in}|m{0.3in}|m{0.3in}|m{0.3in}|m{0.3in}|m{0.3in}|m{0.3in}|m{0.3in}|m{0.3in}|}
\hline
quotient &
~ &
~  &
~  &
~  &
~  &
3  &
~  &
~  &
~  \\\hline
First Remainder &
0 &
4 &
6  &
-20  &
38  &
92  &
-200  &
20  &
~  \\\hline
divisor &
~ &
~  &
~  &
2  &
8 &
0  &
-20  &
~  &
~  \\\hline
\end{xtabular}
%\end{flushright}

\tablehead{}
\hspace*{0.21in}
\begin{xtabular}{|m{1.0in}|m{0.3in}|m{0.3in}|m{0.3in}|m{0.3in}|m{0.3in}|m{0.3in}|m{0.3in}|m{0.3in}|}
\hline
quotient &
~ &
~  &
~  &
~  &
3  &
2  &
~  &
~  \\\hline
Second Remainder &
0 &
-10 &
-20 &
78  &
92  &
-200 &
20  &
~  \\\hline
divisor &
~ &
~  &
2  &
8  &
0 &
-20  &
~  &
~ \\\hline
\end{xtabular}

%\begin{flushright}
\tablehead{}
\hspace*{0.68in}
\begin{xtabular}{|m{1.0in}|m{0.3in}|m{0.3in}|m{0.3in}|m{0.3in}|m{0.3in}|m{0.3in}|m{0.3in}|m{0.3in}|}
\hline
quotient &
~ &
~  &
~  &
3  &
2  &
-5  &
~   \\\hline
Third Remainder &
0 &
20 &
78 &
-8  &
-200  &
20 &
~    \\\hline
divisor &
~  &
2  &
8  &
0 &
-20  &
~  &
~  \\\hline
\end{xtabular}
%\end{flushright}

%\vspace*{-.3in}

%\begin{flushright}
\tablehead{}
\hspace*{1.15in}
\begin{xtabular}{|m{1.0in}|m{0.3in}|m{0.3in}|m{0.3in}|m{0.3in}|m{0.3in}|m{0.3in}|m{0.3in}|}
\hline
quotient &
~ &
~  &
3  &
2  &
-5  &
10    \\\hline
Fourth Remainder &
0 &
-2 &
8 &
0  &
20  &
~    \\\hline
divisor &
2  &
8  &
0  &
-20 &
~  &
~  \\\hline
\end{xtabular}
%\end{flushright}

\vspace*{0.25in}

	As-Samaw'al goes on to deal with the extraction square roots according to a method which he qualifies as being general, because it applies to additive polynomials as well as to those which contain terms to be subtracted\footnote[40]{\textit{al-B\=ahir}, \textit{op. cit.}, text, p.68-70.}. The author attributes its discovery to himself, whereas it is inspired by that of Al-Karaj\={\i}.
	It is at the end of this chapter that As-Samaw'al gives the ``sign rule'', probably in the most complete form that it can be found in the mathematics of Islam\footnote[41]{\textit{Ibid.}, p. 70-71.}. There he mentions not only the rule for multiplying algebraic quantities, but also the rule for their division, and especially their addition and subtraction. But these few lines from As-Samaw'al 's writings again raise the same question regarding the relation between Arab and Indian mathematics. These rules have been enunciated in a manner quite as complete by Bh\=askara II\footnote[42]{Colebrook, \textit{op. cit.}, p. 132-135.}. These two authors are entirely contemporary -- As-Samaw'al may have died about ten years before the Indian author. It is obviously futile to wonder which of them was the inspiration for the other. The question is all the more futile since this ``sign rule'' already appears in the\textit{ Ku\d{t}\d{t}aka-ga\d{n}ita} of Brahmagupta (VIIth century)\footnote[43]{\textit{Ibid.}, p.339.}.
	
	Here we must point out that two translations -- the English rendering by Colebrooke, and the French rendering by Ahmed and Rashed -- have unfortunately translated the Arabic word \textit{n\={a}qis} and the Sanskrit words \textit{rina} or \textit{kshaya} (debt or loss, Colebrooke indicates in a footnote), by the word ``negative''. This term, just like the term zero, carries a number of connotations for us which would have been totally foreign to the Arabs and the Indians. Nothing in their mathematical literature would suggest that they admitted the existence of numbers smaller than zero, especially when one remembers that even the word \textit{sifr} or void did not have for them the meaning we now attach to it, as we have said before. Moreover, the ordering of the ``deficient'' numbers was the reverse of what we would now establish among the negative numbers. For us, the number (-1) is greater than (-5). For Al-Karaj\={\i}, the latter was bigger than the former, since when five objects are \textit{missing} in a sum, more are missing than when only one is missing. It is the same even if we express the idea in terms of debt\footnote[44]{This question has been treated in greater detail in our forthcoming book: \textit{La rationalisation de l'alg\`ebre en pays d'Islam; une nouvelle lecture de l'alg\`ebre arabe.}}.
	The two examples cited above, explicitly attributed to the Indians, as well as others, calculated by the same methods, show that the methods of treating algebraic polynomials among the Arabs was very definitely borrowed from the Indian methods of treating numerical polynomials. Al-Karaj\={\i} does not fail to say so explicitly: ``This method is that which is used for the extraction of known quantities in the Indian reckoning and in other methods of reckoning''\footnote[45]{\textit{al-B\=ahir}, \textit{op. cit.}, text, p. 51.}.
	What about the solution of quadratic equations? In this domain, we do not find any reference to India in the Arabic mathematical literature. Rather it is on the side of Babylon that we must for the origin of the theory of these equations, in view of the obvious similarity between the Arab methods and the Babylonian methods of solving these equations. In any case, it seems that the ``sign rule" of which we have spoken of earlier, and of which Al-Khwarizm\={\i} made mention in the beginning of his work on algebra\footnote[46]{\textit{Kit\=ab al-jabr}, Moucharrafa and Ahmad, eds., Dar al-katib al-`arabi, Cairo, 1968, p. 27-28.} -- at least in the case of multiplication -- must be considered to be of Indian origin.

\section{Summary and Conclusions}
Since the beginning of the history of mathematics in the lands of Islam, whether it be arithmetic or algebra, we immediately see the difficulties which arise when we try to isolate the influence of a single country -- India in this case -- on this mathematics. The study of such an influence transcends the geographical boundaries of the region where it acts, and stretches out to China in the east, and to Greece and Babylon in the west and the north. It is necessary even to include the mathematicians of Alexandria, especially Heron and Diophantus\footnote[47]{Some pages of Brahmagupta's, Al-Khwarizm\={\i}'s and Bh\={a}skara's writings on mensuration are reminiscent of the work of Hero of Alexandria; others, by Brahmagupta, are similar to the writings of Diophantus.}. This study would belong to a field of research not yet formed: that of the comparative history of mathematics. Here we do not have the space -- nor the competence ! -- to write the first chapter of such a study.
	The examination of the connections between Arab and Indian mathematics, even when limited to a few particular examples, has helped us to detect the existence of what could be called centres of mathematical synthesis, melting pots situated at the confluence of various influences. India is indisputably one of them. In its turn it contributed to the formation of another centre of synthesis: that of Arab mathematics from the IXth century onwards. If its influence on the theory of equations -- other than the ``sign rule'' -- can be debated, it is manifest in the domains of arithmetic and of operations on algebraic polynomials. It was finally in the Near and Middle East that the larger part of Greek and Oriental mathematics would be resumed, developed greatly and cast in the definitive form it would take in order to pass into the Occident from the XIIIth century onwards.

\section{Translator's Acknowledgements}
I am deeply indebted to my my three (Indian) French teachers: Ms. Piroj Kirodian (now sadly deceased), Ms. Komala Devarajan -- both of Little Angels' School, Mumbai; and Dr. Raji Ramani, Ruia College, Mumbai, for the instruction I have received from them. Their enthusiasm for the French language and culture was contagious. 

I am grateful to Dominique Tourn\`{e}s, editor of \textit{L'Océan Indien au carrefour des mathématiques arabes, chinoises, européennes et indiennes}, for permission to translate this paper from his collection. I thank Prof. Agathe Keller, for her encouragement, and her help in communicating with the editor and publisher.

\end{document}